# A new elastoplastic-damage model with the correction of stress triaxiality and Lode angle

AbdelkhalAk El Hami[1], Bouchaib Radi[2], David Bassir[3,4] *

[1] *LMN, INSA de Rouen,* 685 Av. de l'Université, 76800 Saint-Étienne-du-Rouvray, France
[2] *LIMII, FST Settat, Morocco*
[3] *CNRS/UMR5060, UTBM, Université de Bourgogne Franche Comte, Rue de Leupe, 90400 Sevenans, France*
[4] *Centre Borelli, Université Paris Saclay, 4 Av. des Sciences, 91190 Gif-sur-Yvette, France*

*\* Corresponding author*: david.bassir@utbm.fr ; david.bassir@ens-paris-saclay.fr

**Abstract:** The classic elastoplastic-damage constitutive model neglects the effects of loading histories. But in fact, more and more experiments results show that the states of stress can significantly affect the response of metals not only in the plasticity but also in the damage evolution. This paper presents an elasto-plastic fully coupled ductile damage constitutive model which considers the corrections of the metal plasticity and the damage evolution through the states of stress. The stress triaxiality and the Lode angle are frequently used to characterize the state of stress. The correction of metal plasticity is achieved by applying a new flow stress model associated with the stress triaxiality and Lode angle. A stress state parameter which depends on the function of stress triaxiality and Lode angle is also proposed. In order to calibrate the damage evolution, the stress state parameter is introduced into the Chaboche damage model. The experiment results of aluminium 2024-T351 are used to validate the proposed constitutive model.

## 1. Introduction

See the end of the document for further details on references. The prediction of plasticity and ductile damage of metal is an important topic in modern manufacturing industrials. During the metal forming processes, the material must have excellent plasticity to flow in the complex dies and the ductile damage should be avoided, in order to obtain a high quality products. During the metal machining, the material must experience the fracture, in order to remove the wastage using various cutting tools. Over the past several decades, some macroscopic models have been developed to predict the material behaviour of plasticity and ductile damage [1, 2, 3, 4, 5, 6]. Although these models have been applied to predict the material behaviour, they still cannot adapt to complex loading histories which obviously exits in the examples mentioned above. Thus, some new constitutive model, which depends on the stress states (tension, shear or compression), is necessary to developed in order to accurately predict the material behaviour in complex loading histories.

The three stress invariants are frequently used to characterize the states of stress [7, 8, 9]. Recent researches on effects of loading histories mainly focused on the Lode angle and stress triaxiality, which can be defined using the three stress invariants. Johnson and Cook [5] considered the stress triaxiality as a part of the criterion of the plastic strain to fracture. Bardet [10] introduced the stress invariants into the constitutive equations to predict the soil and rock behaviour. Zhang [11] took into account the common triaxiality parameter and the Lode parameter in the three-dimensional numerical analyses of a spherical void. Peric [12] applied a pore pressure equation to the three-invariant Cam clay model. Two sets of loading histories are presented: proportional and circular stress paths. Bai [9] presented a general form of asymmetric metal plasticity which considers both the pressure sensitivity and the Lode dependence. The model also validated by experimental results using the specimens of different shapes. Beese [13] developed a phenomenological macroscopic plasticity model for steels. The model is a stress-state dependent transformation kinetics law that accounts for both the effects of the stress triaxiality and the Lode angle on the rate of transformation. The above investigations point out that the stress triaxiality and Lode angle have the ability to characteristic the states of stress and they have the significant effects on material behaviours. However, the fully coupled ductile damage model based on Continuum Damage mechanics (CDM) has never improved to consider the dependence of stress triaxiality and Lode angle, which aims to calibrate the plasticity and damage evolution. The objective of this paper is to develop such a constitutive model in order to take into account the effects of stress triaxiality and Lode angle into the material behaviour of elastoplastic fully coupled ductile damage.

The first aspect, an elastoplastic fully coupled ductile damage constitutive model is built based on the Continuum Damage Mechanics. The damage in CDM theory can be contributed of the nucleation and growth of micro defects (voids and cracks), and their coalescence into macrocracks [14, 15]. Kachanov [16] is the pioneer to characterized this ductile damage by a certain scalar to define the effective stress. Lamaitre and Choboche [4, 6, 17, 18] introduced a damage evolution model using the effective definition of state variables. In this paper, the ductile damage is also assumed to be one of the internal state variables which relates to material behaviour induced by the irreversible deterioration of microstructure.

The second aspect, a new flow stress model is used in this elastoplastic-damage constitutive model in order to calibrate the metal plasticity. This flow stress model considers the effects of stress triaxiality (associated with the size of yield surface) and Lode angle (associated with the shape of yield surface). The parameters of calibration are found using the experiment results in literature [9].

The third aspect, a stress state parameter is proposed, which shows the function of stress triaxiality and Lode angle. It is integrated into the Chaboche damage evolution model and the elastoplastic model. In order to determine the accuracy of the stress state parameter, some experiment results of literatures [19, 7, 20], are used. These experiments are conducted in some special designed specimens who can generate different state of stress. Thus, the stress state parameter can correct the damage evolution in complex loading histories.

The material behaviour of 2024-T351 aluminium is predicted using the proposed constitutive model in this paper. The validation of this constitutive model is implemented in the tensile, shear and compressive tests. The simulation results show that the proposed elastoplastic-damage model with the correction of stress triaxiality and Lode angle can accurately predict the material behaviour under complex loading histories.

## 2. Stress state-dependent elasto-plastic-damage model

This section describes the fundamental formulations for a fully coupled elasto-plastic-damage constitutive model that depends on stress states for ductile materials. The stress state of an isotropic material is firstly defined by applied the parameters of stress triaxiality and Lode angle. Next, the plasticity is corrected by stress state through introducing the stress triaxiality and Lode angle into the von Mises yield criterion. As the same time, the effective state variables considering the ductile damage is modified in order to correct the ductile damage evolution in multi-axial conditions. In this research, the material assumes to be incompressible and has the behaviour of isotropic hardening and ductile damage.

### 2.1. Stress state: stress triaxiality and Lode angle

The stress state at a point inside an isotropic material can be defined by the second order Cauchy stress tensor $\underline{\sigma}$ in all directions. However, the stress tensor $\underline{\sigma}$ changes according to the system of coordinates. Thus, the stress state is frequently described by three stress invariants. The first principal invariant of the Cauchy stress ($\underline{\sigma}$), and the second and third principal invariants ($J_2, J_3$) of the deviatoric part ($\underline{S}_{ij} = \underline{\sigma}_{ij} - \sigma_m \delta_{ij}$) of the Cauchy stress are defined as:

$$\begin{cases} I_1 = \sigma_{11} + \sigma_{22} + \sigma_{33} \\ J_2 = \frac{1}{2} \underline{S}_{ij} : \underline{S}_{ij} \\ J_3 = \det\left[\underline{S}_{ij}\right] \end{cases} \tag{1}$$

where $\delta_{ij}$ is the Kronecker delta. Then, a related set of quantities, ($\sigma_m, \sigma_{eq}, r$), can be defined according to these invariants as:

$$\sigma_m = \frac{1}{3}I_1; \quad \sigma_{eq} = \sqrt{3J_2}; \quad r = \left(\frac{27}{2}J_3\right)^{\frac{1}{3}} \qquad (2)$$

where $\sigma_{eq}$ is equivalent stress. Then, the stress triaxiality, an important stress state parameter [21, 22, 23], can be defined:

$$\eta = \frac{\sigma_m}{\sigma_{eq}} \qquad (3)$$

The Lode angle $\theta$, another important stress state parameter, is related to the normalized third deviatoric stress invariant $\chi$, which is given by:

$$\chi = \left(\frac{r}{\sigma_{eq}}\right)^3 = \cos(3\theta) \qquad (4)$$

and:

$$\theta = \frac{1}{3}\arccos \chi \qquad (5)$$

The normalized third deviatoric stress invariant $\chi$ has the range of $-1 \le \chi \le 1$ and the Lode angle $\theta$ has the range: $0 \le \theta \le \pi/3$.

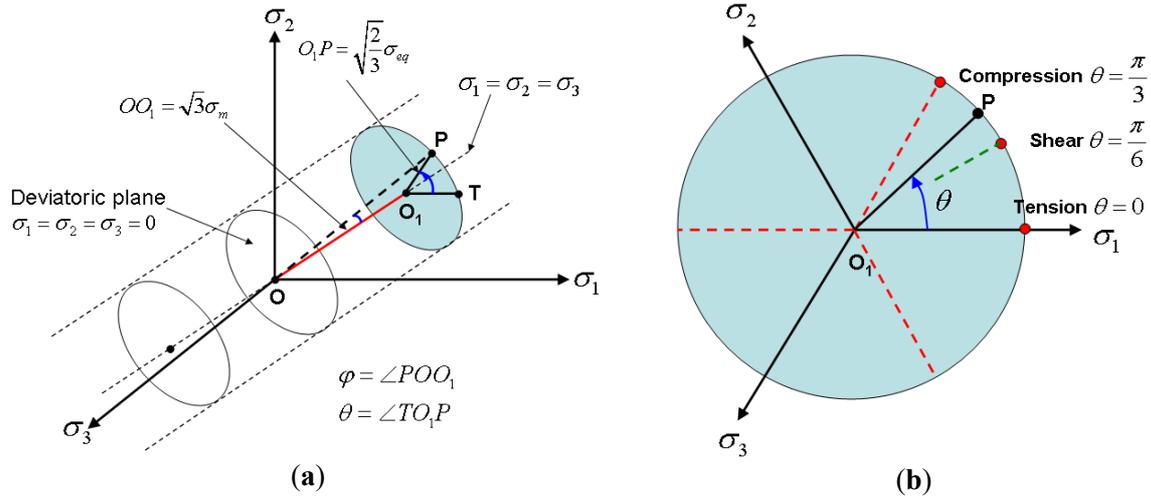

**Figure 1.** Geometrical representation of the stress triaxiality and Lode angle: (**a**) Stress triaxiality and Lode angle; (**b**) Lode angle and stress state.

*2.2. Stress state: stress triaxiality and Lode angle*

The material plasticity can be corrected by using a stress state-dependent flow stress model. The pressure-dependent model is firstly presented by Mohr–Coulomb in order to model the plastic flow of geomaterials and other cohesive-frictional materials. A smooth version of the Mohr–Coulomb yield surface is introduced by Drucker–Prager [24]. It is efficient to determine whether a material has failed or undergone plastic yielding for soil or other pressure-dependent materials. Recently, more and more experimental results prove that not only the stress triaxaility but also the Lode angle significantly affects the plastic flow. Bai and Wierzbicki [9] presented a new model of metal plasticity and fracture which is depends on the pressure and Lode parameter. They explained the dependence of yield conditions on Lode angle by comparing the von Mises and Tresca yield conditions in polar coordinate system.

According Bai's works, a new term considering a linear function of stress triaxiality together with a power function of Lode angle parameter are introduced into the flow stress model. The following flow stress is proposed:

$$\sigma_y(\bar{\varepsilon}^p, \eta, \theta_0) = \bar{\sigma}(\bar{\varepsilon}^p)\left[1 - c_\eta(\eta - \eta_0)\right]\left\{c_\theta^s + (c_\theta^{ax} - c_\theta^s)\left[\theta_0^2 - \frac{(\theta_0^2)^{m+1}}{m+1}\right]\right\}\frac{\delta\bar{\sigma}(\bar{\varepsilon}^p)}{\delta\bar{\varepsilon}^p} \quad (7)$$

where $c_\theta^{ax}$ is a parameter defined by:

$$c_\theta^{ax} = \begin{cases} c_\theta^t & \text{for } \theta_0 \geq 0, \\ c_\theta^c & \text{for } \theta_0 < 0. \end{cases} \quad (8)$$

Eq. (7) defines the shape of yield surface which associates with the plastic strain, pressure and Lode angle parameter. It borrows the form from the Bai's model [9] but it is seems more simple to use.

The Tresca yielding is a hexagon inscribed on the von Mises circle. As presented in Fig. 2(a), the yield surfaces described by Tresca, von Mises and corrected von-Mises models are illustrated. The Tresca yield condition considers the shear state but it neglects the effect of middle principal stress. Mathematically speaking, it must to calculate the principal stress and some differential problems exist in the angular points. The von-Mises is closer to the experimental results but it is not accurate when stress state transfers to shear condition. In this paper, the corrected von Mises is proposed to consider the effects of stress triaxiality and Lode angle. Comparing to Tresca and von-Mises yield, the yield stresses are the same in tension while the shear and compressive stress states are well corrected by the set of parameters: ($c_\theta^t = 1, c_\theta^s = 0.92, c_\theta^c = 1.05, m = 6$). The differential problem is also solved thanks to a high order power function $\frac{(\theta_0^2)^{m+1}}{m+1}$. This function makes the yield surface smooth and differentiable near the point $\theta_0 = \pm 1$, which is also introduced in the Bai's work [9]. The smooth operation using different values of parameter $m$ (non-negative) is clearly presented in Fig. 2(b). A relative larger $m$ is needed to accurately adjust the plasticity when stress state changes.

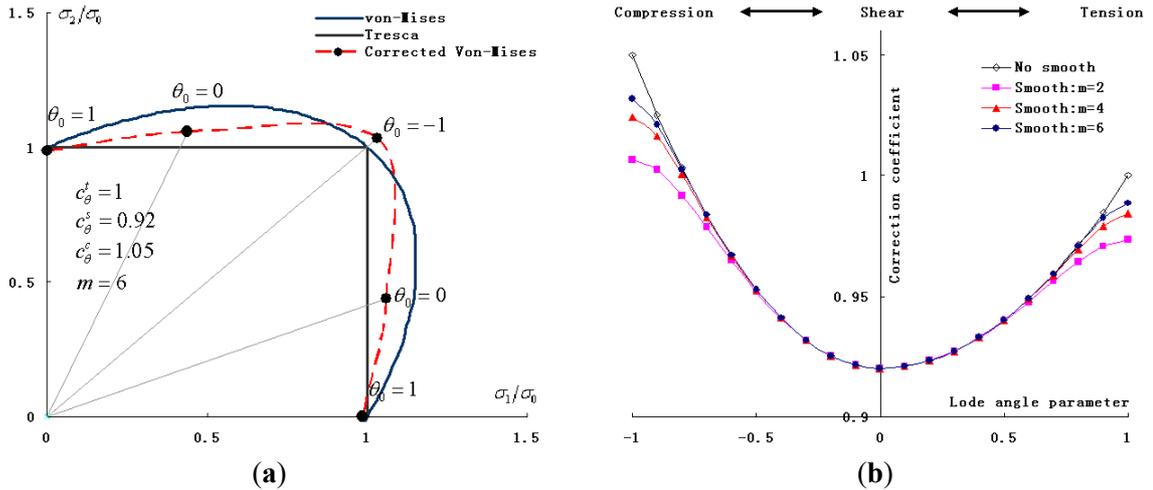

**Figure 2.** Yield surface: von-Mises, Tresca and corrected von-Mises; Calculated coefficients for yield in corrected von-Mises model and the smoothing operation on differential problems using high order power function (Effect of stress triaxiality is neglects: $c_\eta = 0$): (a) Comparison of yield surfaces; (b) Smoothing operation.

*2.3. Effective state variables for multi-axial damage evolution*

As the point mentioned above, the fracture surfaces of ductile material (e.g. metals) have different morphologies when stress state changes. The mechanism of ductile fracture has been proved to be different in tension, shear and compression. The continuum damage mechanism is based on the nucleation, growth, coalescence and linkage of micro-voids, which agrees with the damage evolution in tension. While in the conditions of shear or compression, the damage evolution is closely related to the shear localization between voids or the cracks closures. Hence, the ductile damage model based on the continuum damage mechanics must be corrected by stress state (stress triaxiality and Lode angle parameter) in order to make it available in multi-axial loading conditions.

Unlike the definition of plastic strain to fracture in some models [25, 5, 19, 7, 8, 26, 20, 9], the ductile damage in continuum damage mechanics is related to the nucleation and growth of micro defects (voids and cracks), and their coalescence into macrocracks [1, 2, 3, 4, 6, 27, 28]. It is evaluated by the ductile damage variable D ($0 \leq D < 1$) which is defined by the ratio of effective section area (without the area of microdefects $\tilde{A}$) to whole cross-sectional area $A$, as:

$$D = \frac{A - \tilde{A}}{A} \tag{9}$$

in which: D=0 corresponds that the material is undamaged; D=1 means that the material is fully damaged. In order to couple the effect of the damage variable on the mechanical behaviour, the effective state variables $\left( \tilde{\underline{\sigma}}, \tilde{\underline{\varepsilon}} \right)$ are introduced:

$$\tilde{\underline{\sigma}} = \frac{\underline{\sigma}}{\sqrt{1 - hD}}, \quad and \quad \tilde{\underline{\varepsilon}}^e = \sqrt{1 - hD}\, \underline{\varepsilon}^e \tag{10}$$

where $h$ is the stress state parameter for ductile damage. The classical definition of effective state variables considers that the damage increase is the same in tension, shear and compression. Lemaitre has firstly suggested the use of a regularization variable $h$ to account for crack closure effects and the value of $h$ is often taken close to 0.2 for steels in compression [29]. Generally speaking, the smaller value of stress state parameter gets, the slower damage evolution reaches. In this paper, not only the crack closure effects in compression but also the shear damage mechanism are considered thanks to define the stress state parameter for ductile damage. It is the function of stress triaxiality $\xi$ and Lode angle parameter $\theta_0$, as:

$$h(\eta, \theta_0) = 1 + d_\theta^s + \left[ d_\eta^{ax}(\eta - \eta_0) - d_\theta^s \right] \left[ \theta_0^2 - \frac{\left( \theta_0^2 \right)^{m+1}}{m+1} \right] \tag{11}$$

in which: $\eta_0$ is the reference stress triaxiality and $d_\eta^{ax}(\theta > 0: d_\eta^c$ or $\theta < 0: d_\eta^t)$ is the material constant, which are used to consider the effect of stress triaxiality; $d_\theta^s$ and $m$ are parameters to control the effect of Lode angel parameter. The reference test should be chosen for stress state parameter too, for example the smooth round bar tensile test ($\eta_0 = 1/3$). In this definition, the damage evolution in pure shear condition, controlled by parameter $d_\theta^s$, is assumed to be independent with stress triaxiality. Note that there are three special points respectively correspond to different stress states when the parameter $m$ tends to be infinitely large, as:

$$h(\eta,\theta_0) = \begin{cases} 1+d_\eta^t(\eta-\eta_0) & \text{for} \quad \theta_0 = 1 \\ 1+d_\theta^s & \text{for} \quad \theta_0 = 0 \\ 1+d_\eta^c(\eta-\eta_0) & \text{for} \quad \theta_0 = -1 \end{cases} \qquad (12)$$

In conclusion, Eq. (11) defines the stress state parameter in order to correct damage variable D. It is different to Eq. (7) which is used for the correction of plasticity. Totally, when stress state transfers from tension to compression, the stress state parameter should be reduced in order to prolong the damage evolution while the yield stress must be heightened for strengthening the plastic hardening. The definition of parameter $m$ is the same, which is used to make the yield surface smooth and differentiable near the point $\theta_0 = \pm 1$, as shown in Fig 2(b). So far, the main purpose in this paper must be changed to couple these two equations into the constitutive model in the following sections.

## 2.4. Elastoplastic-damage constitutive equations

The definitions of stress state-dependent yield model and effective state variables based on the stress state parameter make the model powerful in multi-axial loading conditions. To the further implementation, the Eq. (7) and Eq. (11) must be coupled into the elastoplastic-damage constitutive equation. The constitutive relations of state variables are deduced from the elastic free energy, as:

$$\rho\psi_e(\underline{\varepsilon}^e, D) = (1-hD)\left[\mu_e\left(\underline{\varepsilon}^e : \underline{\varepsilon}^e\right) + \frac{\lambda_e}{2}\left(\underline{\varepsilon}^e : \underline{1}\right)^2\right] \qquad (13)$$

where $\rho$ is the material density in the current undamaged configuration, $(\mu_e, \lambda_e)$ are the classical Lame's constants which are calculated with:

$$\lambda_e = \frac{\nu E}{(1+\nu)(1-2\nu)} \qquad \text{and} \qquad \mu_e = \frac{E}{(1+\nu)} \qquad (14)$$

where $E$ and $\nu$ denote the elastic modulus and Poisson's ratio.

According to Eq. (13), the stress-like variables (Cauchy stress tensor and scalar state variable associated D), can be derived by Clausius-Duhem Inequality respectively:

$$\underline{\sigma} = \rho\frac{\partial\psi}{\partial\underline{\varepsilon}^e} = 2\mu_e(1-hD)\underline{\varepsilon}^e + \lambda_e(1-hD)\left(\underline{\varepsilon}^e : \underline{1}\right)\underline{1} \qquad (15)$$

$$Y = -\rho\frac{\partial\psi}{\partial D} = h\left[\mu_e\left(\underline{\varepsilon}^e : \underline{\varepsilon}^e\right) + \frac{\lambda_e}{2}\left(\underline{\varepsilon}^e : \underline{1}\right)^2\right] \qquad (16)$$

where $h$ is the stress state parameter which has been coupled into the Cauchy stress tensor $\underline{\sigma}$ and the damage energy release rate Y. And Eq. (16) shows us that the stress state parameter is a coefficient of Y and it works on changing the damage energy release rate. In other words, the stress state parameter can correct the accumulation speed of damage energy in the different stress states, which is the main conception to define the stress state parameter for multi-axial damage evolution. It should be noted that the damage energy is associated with elastic energy in this model and it can be also extended to be associated with the plastic energy.

The overall strain-rate $\underline{\dot{\varepsilon}}$ in the ductile material can be decomposed into its elastic part $\underline{\dot{\varepsilon}}^e$ and plastic part $\underline{\dot{\varepsilon}}^p$. The elastic part is recoverable while the plastic part is unrecoverable. Hence, the dissipation potential $F$ in an elastoplastic-damage body must be composed of two parts: the yield criterion $f$ and the damage dissipation potential $F_Y$:

$$F = f + F_Y \tag{17}$$

where:

$$F_Y = \frac{\gamma}{(1-hD)^\beta} \frac{1}{\alpha+1} \left(\frac{Y-Y_0}{\gamma}\right)^{\alpha+1} \tag{18}$$

$$f = \frac{\sigma_{eq}}{\sqrt{1-h(\eta,\theta_0)D}} - \sigma_y(\bar{\varepsilon}^p,\eta,\theta_0) \tag{19}$$

Eq. (17) is the Lemaitre damage potential and set (Y0, α, β and γ) are parameters to control the damage evolution. In Eq. (19), $\sigma_{eq} = \sqrt{\frac{3}{2}\underline{S}:\underline{S}}$ is the second invariant of the stress tensor and $\underline{S}$ is the deviatoric part of the Cauchy stress tensor. The stress state-dependent yield stress $\sigma(\bar{\varepsilon}^p,\eta,\theta_0)$ has defined before and is coupled in the yield criterion $f$ directly. To simplify the expression, the stress state parameter and stress state-dependent yield stress are respectively written like $h$ and $\sigma_y$ instead of $h(\eta,\theta_0)$ and $\sigma_y(\bar{\varepsilon}^p,\eta,\theta_0)$ hereinafter.

The evolution of all relationships mentioned above are governed by the plasticity which uses the unique potential dissipation $F$ and the unique yield criterion $f$. According to the rule of normality of the mechanical dissipation, it is possible to define a single plastic multiplier $\dot{\lambda}$ and obtain the formulation:

$$\phi_m = \underline{\sigma}:\dot{\underline{\varepsilon}}^p + Y:\dot{D} - \dot{\lambda}F(\underline{\sigma},Y;D) \tag{20}$$

where the plastic multiplier $\dot{\lambda}$ is a positive scalar and it is determined by consistency condition, as:

$$\dot{f} = \frac{\delta f}{\delta \underline{\sigma}}:\dot{\underline{\sigma}} + \frac{\delta f}{\delta \bar{\varepsilon}^p} \cdot \dot{\bar{\varepsilon}}^p + \frac{\delta f}{\delta D}\dot{D} = 0 \tag{21}$$

The plastic multiplier $\dot{\lambda}$ is deduced:

$$\dot{\lambda} = \frac{1}{h_p}\left\langle \frac{\delta f}{\delta \underline{\sigma}}:\dot{\underline{\sigma}} \right\rangle \tag{22}$$

where $h_p$ is the plastic hardening modulus:

$$h_p = -\frac{\delta f}{\delta \bar{\varepsilon}^p}\dot{\bar{\varepsilon}}^p - \frac{\delta f}{\delta D}\dot{D} \tag{23}$$

The strain associated flux variables $(\dot{\underline{\varepsilon}}^p,\dot{D})$ are derived when the extremums of dissipation potential in Eq. (20) are found:

$$\begin{cases} \dot{\underline{\varepsilon}}^p = \dot{\lambda}\frac{\partial F}{\partial \underline{\sigma}} = \dot{\lambda}\frac{\partial f}{\partial \underline{\sigma}} \\ \dot{D} = \dot{\lambda}\frac{\partial F}{\partial Y} = \dot{\lambda}\frac{1}{(1-hD)^\beta}\left(\frac{Y-Y_0}{\gamma}\right)^\alpha \end{cases} \tag{24}$$

The equivalent plastic strain is updated:

$$\dot{\bar{\varepsilon}}^p = \sqrt{\frac{2}{3}\underline{\dot{\varepsilon}}^p : \underline{\dot{\varepsilon}}^p} = \dot{\lambda}\sqrt{\frac{2}{3}\frac{\partial f}{\partial \underline{\sigma}} : \frac{\partial f}{\partial \underline{\sigma}}} \tag{25}$$

The derivation of differentials with respect to the damage and the plastic strain in Eq. (25) can be expressed by following equations:

$$\frac{\delta f}{\delta D} = \frac{h\sigma_{eq}}{2(1-hD)^{3/2}} \tag{26}$$

$$\frac{\delta f}{\delta \bar{\varepsilon}^p} = -\left[1 - c_\eta(\eta - \eta_0)\right]\left\{c_\theta^s + (c_\theta^{ax} - c_\theta^s)\left[\theta_0^2 - \frac{(\theta_0^2)^{m+1}}{m+1}\right]\right\}\frac{\delta \bar{\sigma}(\bar{\varepsilon}^p)}{\delta \bar{\varepsilon}^p} \tag{27}$$

where $\bar{\sigma}(\bar{\varepsilon}^p)$ is any form of isotropic strain hardening model, for example the power function is used in this paper:

$$\bar{\sigma}(\bar{\varepsilon}^p) = A + B(\bar{\varepsilon}^p)^n \tag{28}$$

The deviation of differential $\frac{\delta f}{\delta \underline{\sigma}}$ is more complex because of the existences of $h$ and $\sigma_y$ which is related to three stress invariants:

$$\frac{\delta f}{\delta \underline{\sigma}} = \frac{\delta f}{\delta \sigma_{eq}}\frac{\delta \sigma_{eq}}{\delta \underline{\sigma}} + \frac{\delta f}{\delta h}\frac{\delta h}{\delta \underline{\sigma}} + \frac{\delta f}{\delta \sigma_y}\frac{\delta \sigma_y}{\delta \underline{\sigma}} \tag{29}$$

The differentials $\frac{\delta f}{\delta \sigma_{eq}}$, $\frac{\delta f}{\delta h}$ and $\frac{\delta f}{\delta \sigma_y}$ in Eq. (29) are easier to obtained. While the Cauchy stress $\underline{\sigma}$ is the second order tensor, so the deviation must aim to each component $\sigma_{ij}$: $\frac{\delta \sigma_{eq}}{\delta \sigma_{ij}}$, $\frac{\delta h}{\delta \sigma_{ij}}$ and $\frac{\delta \sigma_y}{\delta \sigma_{ij}}$

$$\frac{\delta \sigma_{eq}}{\delta \sigma_{ij}} = \frac{3}{2}\frac{s_{ij}}{\sigma_{eq}} \tag{30}$$

$$\frac{\delta h}{\delta \sigma_{ij}} = 2\theta_0\left[1-(\theta_0^2)^m\right]\left[d_\eta^{ax}(\eta-\eta_0) - d_\theta^s\right]\frac{\delta \theta_0}{\delta \sigma_{ij}} + d_\eta^{ax}\left[\theta_0^2 - \frac{(\theta_0^2)^{m+1}}{m+1}\right]\frac{\delta \eta}{\delta \sigma_{ij}} \tag{31}$$

$$\begin{aligned}\frac{\delta \sigma_y}{\delta \sigma_{ij}} &= \bar{\sigma}(\bar{\varepsilon}^p)\left[1-c_\eta(\eta-\eta_0)\right]2\theta_0(c_\theta^{ax}-c_\theta^s)\left[1-(\theta_0^2)^m\right]\frac{\delta \theta_0}{\delta \sigma_{ij}} \\ &\quad -\bar{\sigma}(\bar{\varepsilon}^p)c_\eta\left\{c_\theta^s + (c_\theta^{ax}-c_\theta^s)\left[\theta_0^2 - \frac{(\theta_0^2)^{m+1}}{m+1}\right]\right\}\frac{\delta \eta}{\delta \sigma_{ij}}\end{aligned} \tag{32}$$

which Eq. (30) is the normal direction of flow stress without considering stress state in ductile damage and plasticity. The differentials in Eq. (31) and (32) can be finished by the following equations:

$$\frac{\delta \eta}{\delta \sigma_{ij}} = -\frac{3\eta}{2\sigma_{eq}^2} s_{ij} \tag{33}$$

$$\frac{\delta \theta_0}{\delta \sigma_{ij}} = \frac{9}{\pi \sigma_{eq} \sin 3\theta} \left( \frac{\cos 3\theta}{\sigma_{eq}} s_{ij} - \frac{3}{\sigma_{eq}^2} s_{ik} s_{kj} \right) \tag{34}$$

Although the plasticity is pressure sensitive, experiments show that the plastic dilatancy of metals is negligible [30, 9]. For the assumption of plastic incompressibility, the effects of hydrostatic stress tensor on plastic flow direction in above equations should be removed, which means $\frac{\delta f}{\delta \sigma_{ij}} = \frac{\delta f}{\delta s_{ij}}$ in the proposed constitutive equations. Therefore, what is used in the present paper is a flow rule with deviatoric associativity. Then, the plastic modulus is obtained by substituting the equations from (24) to (34) for the terms in Eq. (23).

Furthermore, the elastoplastic modulus can be derived by calculating the time derivative of stress tensor ($\underline{\sigma} = 2\mu_e (1-hD) \underline{\varepsilon}^e$):

$$\dot{\underline{\sigma}} = \frac{\partial \underline{\sigma}}{\partial \underline{\varepsilon}^e} : \dot{\underline{\varepsilon}}^e + \frac{\partial \underline{\sigma}}{\partial D} \dot{D} \tag{35}$$

where $\dot{\underline{\varepsilon}}^e = \dot{\underline{\varepsilon}} - \dot{\underline{\varepsilon}}^p = \dot{\underline{\varepsilon}} - \dot{\lambda} \frac{\partial f}{\partial \underline{\sigma}}$. Considering the consistency condition in Eq. (20), one can obtain:

$$\dot{\lambda} = \frac{1}{H_{ep}} < \frac{\partial f}{\partial \underline{\sigma}} : \frac{\partial \underline{\sigma}}{\partial \underline{\varepsilon}^e} : \dot{\underline{\varepsilon}}^e > \tag{36}$$

and

$$H_{ep} = h_p + \frac{\partial f}{\partial \underline{\sigma}} : \frac{\partial \underline{\sigma}}{\partial \underline{\varepsilon}^e} : \frac{\partial F}{\partial \underline{\sigma}} - \frac{\partial f}{\partial \underline{\sigma}} : \frac{\partial \underline{\sigma}}{\partial D} \frac{\partial F}{\partial Y} \tag{37}$$

All differentials in Eq. (37) are known and $H_{ep}$ gives out the elastoplastic modulus. Finally, the plastic multiplier $\dot{\lambda}$ is solved and other state variables can be updated.

In conclusion, the constitutive equations developed in this paper is a elastoplastic fully coupled ductile damage model which corrects the plasticity and damage evolution according to stress states in the multi-axial loading conditions. It is based on the continuum damage mechanics but has great improvements. The corrected yield criterion and the modified definition of effective state variables have been fully coupled into the constitutive equations. This is the main difference with the models which define the plastic strain to fracture [19, 8, 9]. It must be pointed out that, the plasticity correction and the damage evolution are two different aspects to affect the constitutive relation in this model. The ductile damage evolves throughout the whole procedure of material deformation, from elastic stage to final fracture; while the material plasticity is only corrected through yield stress in plastic stage. Hence, the damage evolution under multi-axial loading which matches to different damage mechanism and the strain hardening along different loading paths can be studied separately.

### 3. Results and validations

To implement this elastoplastic-damage constitutive model into finite element codes, the equations described in the second section are integrated into Abaqus user subroutine VUMAT. The elastic prediction-return mapping algorithm with an operator splitting methodology is used [31, 32, 27, 18] to calculate the plastic multiplier. The fully implicit Euler method is used to update the internal variables, since it contains the property of absolute stability and the possibility of appending further equations to the existing system of nonlinear equations [31, 33].

To apply this numerical methodology, this section focuses on numerically studying the responses of a representative elementary volume under complex stress states and experimentally validating the proposed model in the various tests of 2024-T351 aluminium alloy [22, 23, 34, 35, 37, 38, 39, 40, 41, 42, 43, 44] The responses of a representative elementary volume are researched in order to numerically study the abilities of proposed constitutive model, especially on the aspects of multi-axial plasticity correction and ductile damage evolution. The material parameters for 2024-T351 aluminium alloy are confirmed, which is composed of two steps: elastoplastic-damage parameters are found in the reference tensile test of round smooth bar; the parameters of stress state correction are found in various shear and compressive tests. The validations focus on the tests of tension, shear and compression using various specially designed specimens in literatures [7, 8, 9, 19, 20].

### 3.1. Material responses

Elastoplastic-damage response is the basic behaviour of ductile metals. In the present constitutive model, the parameters A, B and n is used to describe the isotropic strain hardening for metal plasticity while the parameters Y0, α, β and γ are used to control the damage evolution. Neglecting the effect of stress state, a baseline stress-strain curve can be got using these seven parameters. A single experiment test is needed to determine these parameters and this test will be used as the reference test for the further corrections of plasticity and ductile damage.

The tension test of a smooth round bar is used in this paper as the reference test and its force-displacement curve can be found in the experiment of Bai's work [9]. The experimental true stress-strain curve can be computed using these experimental results, as shown in Fig. 3(a). The material properties, strain hardening and damage evolution parameters used for fitting the true stress-strain curve are shown in Tab. 1. Based on these parameters, the material behaviour in the same stress state is predicted by investigating the response of middle integrate point in tensile test of three hexahedral solid elements. It is clear that the coupled model successfully predicts a necking points ( $\bar{\varepsilon}^p = 0.47, \bar{\sigma}_{\max} = 792\text{MPa}$ ) and a final fracture at $\bar{\varepsilon}^p = 0.55$. It is also clear that the degradation of material stiffness is induced by the evolution of ductile damage evolution in Fig. 3(b). Note that, only the date point before specimen necking is used for fitting the stress-strain curve because the transformation equations are valid only up to necking initiation. And this paper looks upon the necking initiation point as the equivalent plastic strain $\bar{\varepsilon}_f^p$, for example $\bar{\varepsilon}_f^p = 0.47$ in this simulation. Comparing to the results predicted by power law in Abaqus, the fully coupled damage constitutive model is powerful to predict the material behaviours including the elastic, plastic and ductile damage.

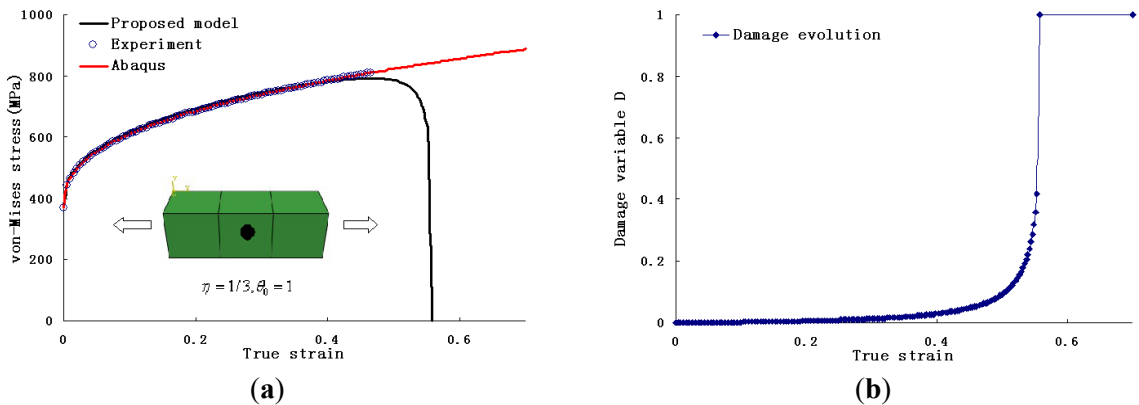

**Figure 3.** Elastoplastic-damage response of 2024-T351 aluminium alloy in the reference test: smooth round bar tensile test ( $\eta = 1/3, \theta_0 = 1$ ): **(a)** True stress-strain curve; **(b)** Damage evolution.

**Table 1.** Elastoplastic-damage parameters and stress state correction constants for 2024-T351 aluminium.

| Material properties | E(MPa) | v | A(MPa) | B(MPa) | n | Y0(MPa) | α | β | γ(MPa) |
|---|---|---|---|---|---|---|---|---|---|
| | 71150 | 0.3 | 370 | 620 | 0.396 | 0 | 2 | 1 | 12.8 |
| Correction constants | $c_\eta$ | $c_\theta^t$ | $c_\theta^s$ | $c_\theta^c$ | $d_\eta^t$ | $d_\theta^s$ | $d_\eta^c$ | m | $\eta_0$ |
| | 0.09 | 1 | 0.855 | 0.9 | 1.3 | 0.55 | 0.6 | 6 | 0.4 |

In above reference test, the stress state is fixed to $\eta = 1/3$ and $\theta_0 = 1$. The next step is to determine the correction constants in proposed material model for correcting the plasticity and ductile damage when stress state changes. After summarizing and analysing the experimental results in the literatures [19, 7, 8, 20, 9] the correction constants used for correcting plasticity and equivalent plastic strain to fracture are found as in Tab. 1. There are four constants ($c_\eta$, $c_\theta^t$, $c_\theta^s$ and $c_\theta^c$) are used to calibrate the plasticity while there are three constants ($d_\eta^t$, $d_\theta^s$ and $d_\eta^c$) are used to calibrate the equivalent plastic strain to fracture. In both two corrections, constant $\eta_0$ is the stress triaxiality in reference test and constant $m$ is the smooth operator which makes the constitutive equations differentiable.

To study the material responses in multi-axial loading conditions, the shapes of two assistant elements in above three hexahedral elements must be changed in order to generate the different stress states in the middle integrate point. For example, in the tension of three hexahedral elements with same shapes (length of elements' edges are the same: L), the stress state in the middle element is $\eta = 1/3$ and $\theta_0 = 1$; while in the tension of three hexahedral elements which change the dimension of loading edges to 1.36L, the stress state in the middle element is $\eta = 0.69$ and $\theta_0 = 1$. A series tests have been designed and conducted based on this method. The material responses with the correction of stress states in tension are shown in Fig. 4. In these tests, the Lode angle parameter is computed to one and the increased stress triaxialities are obtained. One can observe that, the material is harder while the equivalent plastic strain to fracture is higher at low stress triaxiality conditions. These results show us that the trend of correction is in perfect accord with the experimental results. The further validations should be conducted by comparing the plasticity and plastic strain to fracture in various special designed tests.

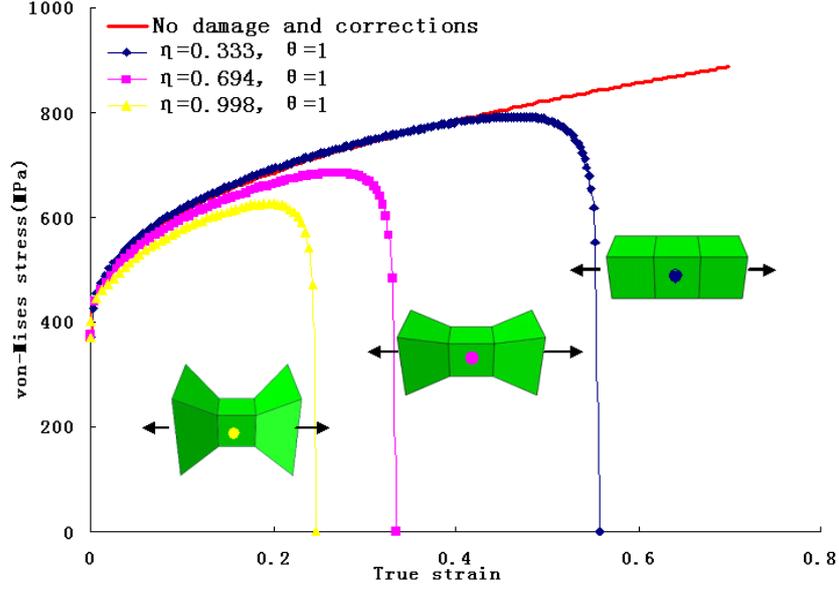

**Figure 4.** Material responses with the correction of stress states in tension

Although these tests are conducted in tension stress state, the stress state parameter and plastic corrector can be derived with the help of above simulation results. For example, the relationship between stress state parameter and equivalent plastic strain to fracture has been reflected in Fig. 4. A power function can be used to describe this relationship, as:

$$\bar{\varepsilon}_f^p = 0.44717 * h^{-1.72555} \tag{38}$$

During the simulation, the stress state parameter is calculated firstly according to different stress triaxialities and Lode angles using the parameters in Tab. 1. Then the damage locus can be predicted through using Eq. (38), as illustrated in Tab. 2. A series experiments has been conducted by Bao [8] and the relationship between equivalent plastic strain to fracture and stress triaxiality has been investigated for 2024-T351 aluminium alloy. When comparing the predicted damage locus to the experimental ones, one can find that the dates agree well with each other, as presented in Fig. 5. Simultaneously, the computing results also have a cut-off value of stress triaxiality which is similar to the experiment results. Besides the stress triaxiality, the lode angle (shear condition) has the significant effects on damage evolution. For example, the fracture generates easily when the Lode angle trends to zero (pure shear) although the stress triaxiality is very negative. That's why it is necessary to develop the fully coupled constitutive model with the correction not only on stress triaxiality but also on Lode angle. The importance of the correction on Lode angle can be investigated in the tension test of flat grooved specimens and in the upsetting test of cylindrical specimens.

**Table 2.** Damage locus predicted by proposed constitutive model in different stress states generated by different specimens.

| Specimen description | $\eta$ | $\theta$ | $h$ | $\bar{\varepsilon}_f^p$ |
|---|---|---|---|---|
| Smooth round bar, tension | 0.4014 | 0.9992 | 1.004109 | 0.463876 |
| Round large notched bar, tension | 0.6264 | 0.9992 | 1.327638 | 0.274216 |
| Round small notched bar, tension | 0.9274 | 0.9984 | 1.666339 | 0.185272 |
| Simple shear | 0.0124 | 0.0355 | 1.54867 | 0.210227 |
| Combination of shear and tension | 0.1173 | 0.3381 | 1.44491 | 0.236951 |
| Pipe, tension | 0.3557 | 0.9186 | 1.035957 | 0.420727 |
| Cylinder($d_0/h_0$=0.5), compression | -0.2780 | -0.7715 | 0.980001 | 0.450136 |
| Cylinder($d_0/h_0$=0.8), compression | -0.2339 | -0.6999 | 1.093852 | 0.375636 |
| Cylinder($d_0/h_0$=1.0), compression | -0.2326 | -0.6794 | 1.120542 | 0.360906 |
| Cylinder($d_0/h_0$=1.5), compression | -0.2235 | -0.6521 | 1.156684 | 0.342403 |
| Round notched, compression | -0.2476 | -0.8941 | 0.79903 | 0.630017 |

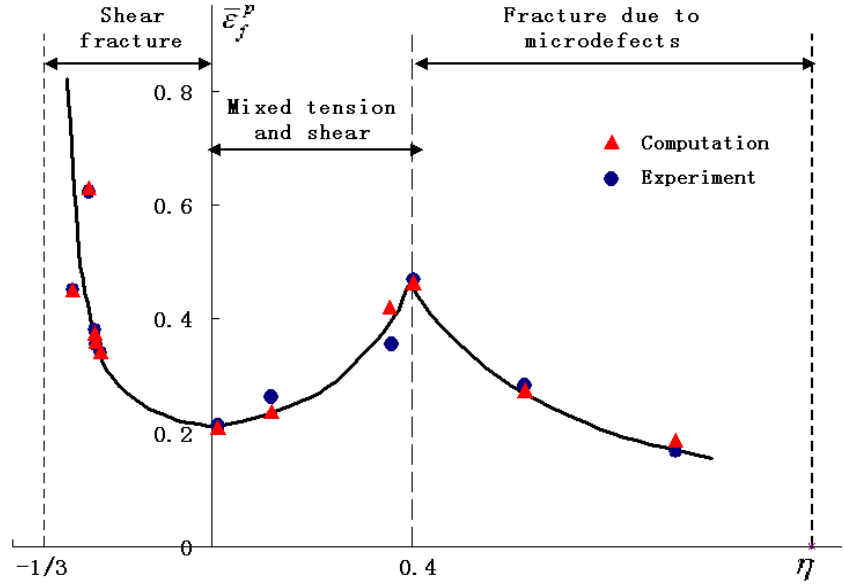

**Figure 5.** Comparison of damage locus: computing results and experimental results [8, 20, 9].

In conclusion, the material responses in complex loading conditions are well predicted through using this fully coupled elastoplastic-damage model with the correction of stress triaxiality and Lode angle. The plasticity and ductile damage evolution for 2024-T351 aluminium alloy have been well corrected in tension, shear and compression. In particular, the plastic strain to fracture used in this paper is the necking point and the stress-strain curves after necking have not been mentioned. There are two aspects are considered: the absence of the force-displacement curve after necking and the period from necking to fracture is very short for metals. For the further researches in some deep necking material, the proposed constitutive model has the ability to predict the degradation of material stiffness, especially continuously predict the residual stiffness after necking (the complex stress state transformation and the controllable damage evolution). In this paper, the behaviour of 2024-T351 aluminium alloy is predicted and the validations will be conducted for different specimens in tension, shear and compression in the next section.

*3.2. Validations in tension, shear and compression*

Generally, the calibrations of correction parameters for plasticity and ductile damage evolution are complicated works, which need a large number of experimental dates. As the work mentioned above, the smooth round bar tensile test ($\eta \approx 1/3$ and $\theta_0 \approx 1$) is needed firstly to determine the basic stress-strain curve (coupled elastoplastic and damage). According to this reference test, the notched round bars tensile tests (only changes the stress triaxiality) should be used to calibrated the pressure effect on plasticity $c_\eta$ and ductile damage $d_\eta^t$ without considering the Lode angle. The third test is the tensile test of flat grooved plates ($\eta \approx 0.66$ and $\theta_0 \approx 0$) or simple shear ($\eta \approx 0$ and $\theta_0 \approx 0$) which can be used to calibrate the parameter $d_\theta^s$ in damage evolution and $c_\theta^s$ in plasticity. The fourth step is to calibrate the parameters $d_\eta^c$ in damage evolution and $c_\theta^c$ in plasticity using the cylindrical specimens' upsetting test (negative stress triaxiality and Lode angle parameter). On the basis of previous Bai's and Bao's works, these parameters have been calibrated by considering the responses of a representative elementary volume in the last section. This section, emphasis is put on the validation of proposed constitutive model using the calibrated parameters in Tab. 1.

Validations in tension are conducted in the tensile simulation of a smooth round bar and a notched round bar. The dimensions and finite element discretization are presented in Fig. 6(a). In these tests, the Lode angle parameter keeps to one while the stress triaxiality increases following with the decrease of notch radius. The smooth round bar tension is the reference test and the constitutive model is validated in notched round bar tension which has different stress triaxialities. Validation in shear is conducted in the tensile simulation of a flat grooved plate whose dimension is shown in Fig. 6(b). For symmetry condition, 1/8 FE model of full plate is used to save the computing scale. In this test, the stress triaxiality in the notched region is similar to notched round bar but the Lode angle parameter is near to 0. The correction with Lode angle parameter can be validated through this test. Validation in compression focus on the cylindrical specimen upsetting test, as illustrated in Fig. 6(c). The stress triaxiality effects on plasticity and ductile damage are validated in tensile tests, which are also inherited in its negative condition. The main purpose of upsetting test is to validate the stress state condition when Lode angle is near to -1.

All of above specimens are discretized by 8 node hexahedral element. The finite element models are solved with the help of Abaqus/explicit solver and material user subroutine Vumat. The minimum mesh sizes used in these four specimens (Fig. 6 from left to right) are 0.5mm, 0.1mm, 0.1mm and 0.15mm respectively. To diminish a possible effect of element size, the damage evolution is configured to associate with the element characteristic length. Actually, the accumulation rate of elastic energy is increased following with the decrease of element characteristic length. A linear relationship between element size and damage evolution is used to guarantee the same damage accumulation in different element size from 0.5mm to 0.1mm. For neglecting the strain rate effect, the same target time increment for mass scaling, is applied in these simulations in order to save computing times. Note that, the friction coefficient in upsetting test is fixed to 0.1 which is also used in literature [9].

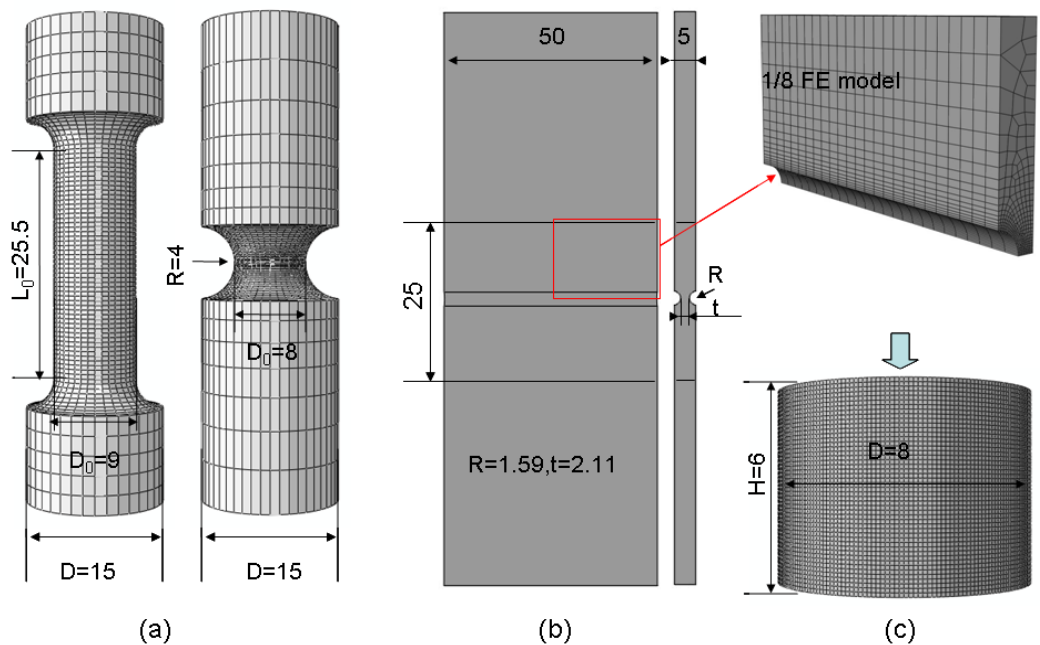

(a)          (b)          (c)

**Figure 6.** Dimensions and finite element discretization specimens ：(a) Smooth and notched round bar tension; (b) Flat grooved plate tension and (c) Cylindrical specimen upsetting.

The validations of proposed model should be concentrated in comparing the force-displacement curves between numerical simulations and experimental results. The comparison of force-displacement curves between experimental results and simulation results in smoothed round bar tensile test is shown in Fig. 7. These curves show us that the simulation results agree well with the experiment ones. Since this test is regarded as the reference test, the plasticity and damage evolution have no significant corrections in the curve predicted by proposed constitutive model. The force-displacement curve predicted by Bai [9] for the same test is also presented in this figure. One can observe that the ductile damage induced degradation of material stiffness is clearly predicted by the proposed model in this paper, even after the necking point. It is proved that the fully coupled ductile damage model is more suitable to simulate the whole material behaviour in the tensile test. The comparison of force-displacement curves between experimental results and simulation results in notched round bar tension is shown in Fig. 8. Without any corrections, one can see the differences between numerical simulations and experiments in force responses. In this figure, the proposed model performs well in predicting the material behaviour in this test. Indeed, the stress triaxiality correction plays important role in correcting not only the plasticity but also the fracture locus in the notched round bar tension.

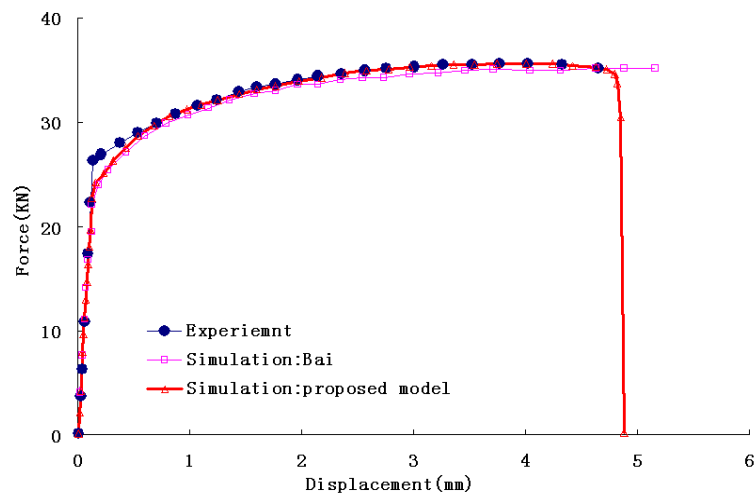

**Figure 7.** Comparison of force-displacement curves between experimental results and simulation results in smoothed round bar tensile test.

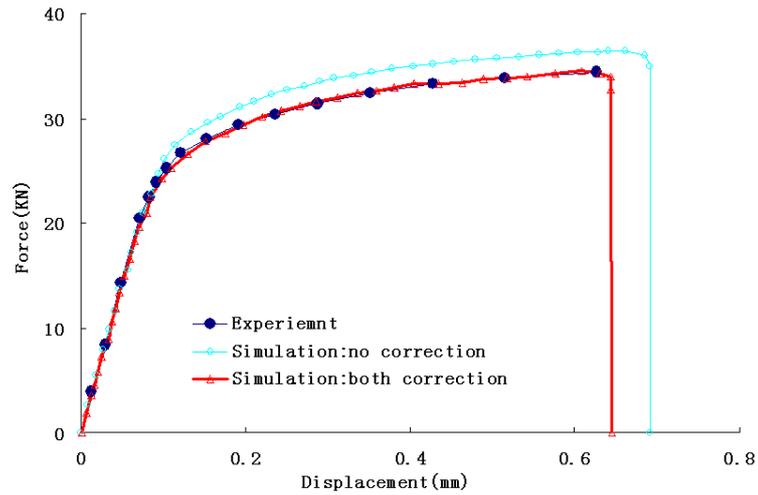

**Figure 8.** Comparison of force-displacement curves between experimental results and simulation results in notched round bar tensile test.

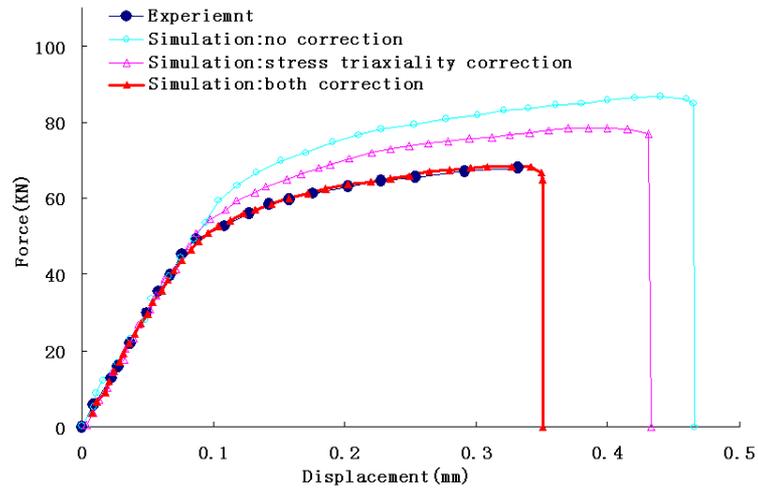

**Figure 9.** Comparison of force-displacement curves between experimental results and simulation results in flat grooved plate tensile test.

The comparison of force-displacement curves between experimental results and simulation results in flat grooved plate tensile test is shown in Fig. 9. The predicted curve without corrections is very different with experiment results in the responses of plasticity and fracture locus. The predicted curve using stress triaxiality correction approaches to experimental results but also has almost 16% errors in plasticity and 23% errors in fracture locus from this figure. That's to say the correction in plasticity and fracture locus according to stress triaxiality is not enough to correct the material behaviour in shear condition for this kind of aluminium alloys. When the Lode angle parameter effect is taken into account, the force-displacement curve using both corrections agrees well with the experiment results.

The comparison of force-displacement curves between experimental results and simulation results in cylindrical specimen upsetting test is shown in Fig. 10. The results show us that the predicted curve without stress state correction has the deviation after the stage (a) in Fig. 10 when the material starts to plastic deformation. These errors have been dispelled through using the both correction constitutive model. It is useful to analyze the damage behaviour in this upsetting test. Through the Fig. 10, the damage is found firstly at the surface center of cylindrical specimen in the stage (b). The stress triaxiality changes between -0.3 and -0.2 and the lode angle stays around -0.7 before this stage. Following with the accumulation of damage, the lode angle at this point increases to zero slowly from stage (c) to (d) and two cross shear band are formed. The microcracks then generate in the surface of specimen and some parallel cracks form until to stage (e). This example proves that the shear induced damage can be well simulated in the upsetting test using proposed model.

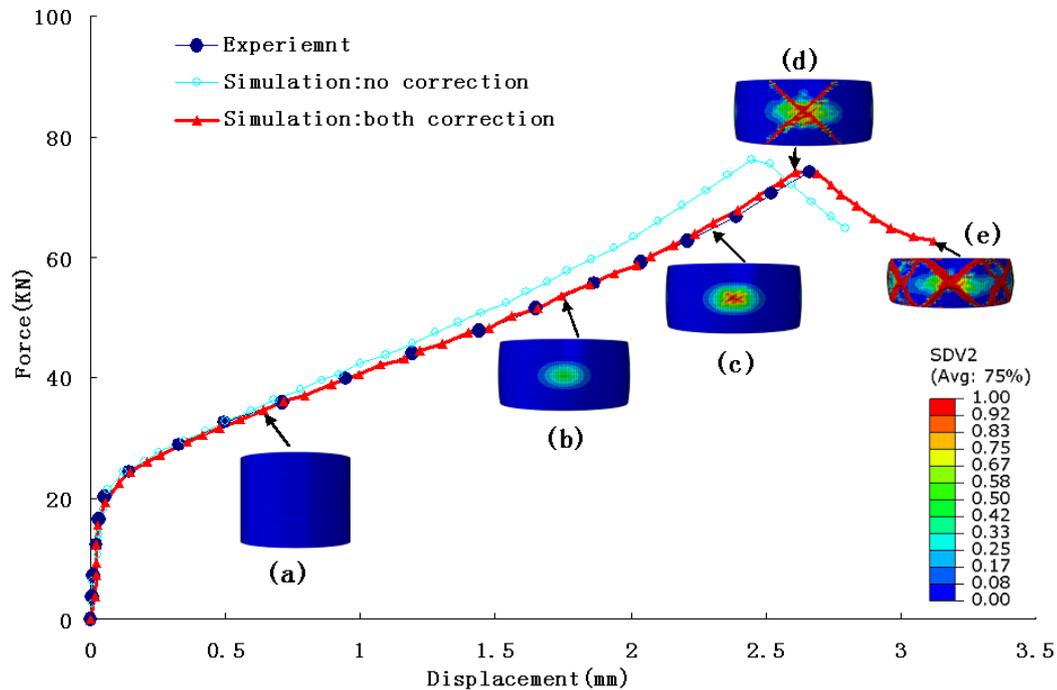

**Figure 10.** Comparison of force-displacement curves between experimental results and simulation results in cylindrical specimen upsetting test.

To conclusion, the proposed elastoplastic-damage model corrected by stress triaxiality and Lode angle is proved to be effective in predicting the material behaviour under complex loading histories. And it is necessary to calibrate the material behaviour in multi axial loading conditions for 2024-T351 aluminium before this alloy is chosen and applied in some metal forming processes.

## 4. Conclusion

A fully coupled elastoplastic and damage constitutive model is developed base on the continuum damage mechanics in this paper. Through correcting the yield stress and effective definition with damage, this model has the abilities to predict material behaviour under complex loading conditions.

The 2024-T351 aluminium alloy is investigated in this paper and the experimental results in literatures prove that its plasticity and fracture locus are sensitive to stress tri-axiality and Lode angle. With the help of the proposed constitutive model proposed, the material behaviour of 2024-T351 aluminium alloy is well simulated in different stress state conditions. The validations are conducted in smooth/notched round bare tensile tests, flat grooved plate tensile test and cylindrical specimen upsetting test. In the future works, this constitutive model will be used to other materials and applied to some metal forming simulations.